\documentclass[11pt,twoside]{article}
\usepackage{graphics,latexsym,amssymb,amsmath}
\RequirePackage{graphics,epsfig,psfrag}
\usepackage{hyperref}
\hypersetup{
    bookmarks=true,         
    unicode=false,          
    pdftoolbar=true,        
    pdfmenubar=false,        
    pdffitwindow=false,     
    pdfstartview={FitH},    
    pdftitle={Every knot is a billiard knot},    
    pdfauthor={PVK, DP},     
    pdfsubject={MEGA11},   
    pdfcreator={PVK},   
    pdfproducer={Producer}, 
    pdfkeywords={keyword1} {key2} {key3}, 
    pdfnewwindow=true,      
    colorlinks=true,       
    linkcolor=black,          
    citecolor=black,        
    filecolor=black,      
    urlcolor=black           
}
\def\QQ{{\bf Q}} 

\def\pn{\medskip\par\noindent}
\def\Frac#1#2{{\displaystyle{{#1}\over{#2}}}}

\def\[#1\]{\begin{eqnarray}#1\end{eqnarray}}
\def\$#1\${\begin{eqnarray}#1\end{eqnarray}}

\def\phi{\varphi}

\def\mod{{\rm mod}\ }

\newcommand{\Pf}{{\em Proof}. }

\newcommand{\EPf}{\hbox{}\hfill$\Box$\vspace{.5cm}}

\def\Frac#1#2{{\displaystyle{\frac{#1}{#2}}}}
\def\phi{\varphi}

\newtheorem{theorem}{Theorem}

\newtheorem{remark}[theorem]{Remark}
\newtheorem{lemma}[theorem]{Lemma}
\newtheorem{proposition}[theorem]{Proposition}
\newtheorem{corollary}[theorem]{Corollary}
\newtheorem{thm*}[theorem]{Theorem}
\date{\today}
\begin{document}
\pagestyle{myheadings}
\markboth{ D. Pecker, \today}{{\em Poncelet's theorem and billiard knots}}
\title{Poncelet's theorem and Billiard knots}
\author{ Daniel Pecker}
\maketitle
\begin{abstract}
Let  $D$ be  any elliptic right cylinder. We prove that  every type of knot can be realized
as the trajectory of a ball in $D.$
This proves a conjecture of Lamm and gives a new proof of a conjecture of Jones and Przytycki.
We use Jacobi's proof of Poncelet's theorem by means of elliptic functions.
\pn {\bf keywords:} {Poncelet's theorem, Jacobian elliptic functions,  Billiard knots ,
Lissajous knots,  Cylinder  knots}
\pn
\end{abstract}
\begin{center}
\parbox{12cm} {\small
\tableofcontents
}
\end{center}
\vspace{1cm}
\section{Introduction}

The Poncelet closure theorem is one of the most beautiful theorems in geometry. It says
that if there exists a closed polygon inscribed in a conic $ E$ and circumscribed about
another conic, then there exist infinitely such polygons, one with a vertex
at any given point of $E.$

When the conics are concentric circles the proof is very simple, each Poncelet polygon
is obtained by rotating any one of them.
What makes Poncelet's theorem great, is that it is impossible to generalize this simple proof.

Poncelet's proof ( \cite{Po}) uses "pure" projective geometry, see \cite{Be, DB, Sa} for
 modern proofs along these lines.
Shortly after the publication of Poncelet's book, Jacobi ( \cite{Ja} ) gave a proof
by means of Jacobian elliptic functions. He discovered what is now called a uniformization
of the problem by an elliptic curve. Most modern developments and generalizations
follow Jacobi's proof, see  \cite{Lau, GH, BKOR, Sc, LT}  .

A beautiful example of Poncelet polygonal lines is given by elliptic billiards.
  If a segment of a billiard trajectory in an ellipse $E$ does not
 intersect the focal segment $ [F_1 F_2]$ of $E,$  then there exists an
  ellipse $C$ called a caustic, such that the trajectory is a Poncelet polygonal line
  inscribed in $E$ and circumscribed about $C,$ see \cite{St, T, LT}.

On the other hand,  Jones and Przytycki defined billiard knots as
   periodic billiard trajectories without  self-intersections in a three-dimensional billiard.
They proved that billiard knots in a cube are very special knots,
the Lissajous knots. They also conjectured that every knot is a billiard knot in some
convex polyhedron.
(\cite{JP}, see also \cite{La2,C,BHJS,BDHZ, P}).

Lamm and Obermeyer \cite{La1,LO} proved that
 not all knots are billiard knots in a cylinder. Then Lamm conjectured that there exists
an elliptic cylinder containing all knots as billiard knots (\cite{La1,O}).
It is easy to see that Lamm's conjecture implies the conjecture of Jones and Przytycki:
if $K$ is a billiard knot in a convex set, then it is also a billiard knot in the polyhedron
delimited by the tangent planes.
Dehornoy constructed in \cite{D} (see also \cite{O})  a billiard which contains all knots,
 but this billiard is not convex.
\pn
In this paper we will use Jacobi's method to study billiard trajectories in a right cylinder
with an elliptic basis. We obtain a proof of Lamm's conjecture.
Our result is more precise:

\pn
{\bf Theorem \ref{th:billiardkl}}
{\em Let $E$ be an ellipse which is not a circle,
 and let $D$ be the elliptic cylinder $ D= E \times [0,1].$
Every knot (or link) is a billiard knot (or link)  in $D$.}
\pn

 Billiard trajectories in an ellipse are introduced in section 2.
We show that an elementary theorem of Poncelet implies the existence of a caustic.
We also show that Poncelet polygons in a pair of nested ellipses are projections of torus knots.
Then, by a theorem of Manturov (\cite{M}), we deduce that every knot has a projection which is
a billiard trajectory in an ellipse.

In section 3, we recall the basic definitions and properties of the Jacobian elliptic functions
 sn$(z) $ and cn$(z).$ Then we give the Hermite--Laurent version of Jacobi's proof, which is
 based on Jacobi's uniformization lemma.

In section 4, we use Jacobi's lemma to compute the coordinates of the crossings and vertices of
a billiard trajectory in $E.$
We deduce that if the number of sides of a periodic billiard trajectory is odd, then it is
generally completely irregular.
This means that if we start at any vertex, $1$ and the arc lengths of the other vertices
and  crossings
  are linearly independent over $ \QQ.$
We also see how our proofs generalize in the link case.

In section 5,
we use  Kronecker's density theorem to obtain  our main result.
The same strategy was used in \cite{KP2} to give an elementary proof of
 the Jones--Przytycki conjecture.
 There is
  another application of Kronecker's theorem to the construction of knots in  \cite{KP1}.

\section{Billiard trajectories in an ellipse}\label{star}

The study of billiard trajectories in an ellipse was introduced by Birkhoff
 in 1927 (\cite{Bi}), see also \cite{T} for a modern   exposition  of the subject.

\subsection{Some elementary facts}
The following elementary theorem is due to Poncelet (\cite{Po,Be}).

\begin{theorem} {\bf (The second little Poncelet theorem)}

Let $E$ be an ellipse (or a hyperbola) with foci $F_1$ and $F_2.$
Let $PM_1$ and $PM_2$ be the tangents to $E$ at the points $M_1$ and $M_2.$
Then  the angles $\widehat{M_1PM_2}$ and $ \widehat{F_1PF_2}$ have the same bisectors.
\end{theorem}
\psfrag{P}{$P$}
\psfrag{F1}{$F_1$}
\psfrag{F2}{$F_2$}
\psfrag{M1}{$M_1$}
\psfrag{M2}{$M_2$}
\psfrag{Fp1}{$F'_1$}
\psfrag{Fp2}{$F'_2$}
\begin{figure}[th]
\begin{center}
\epsfig{file=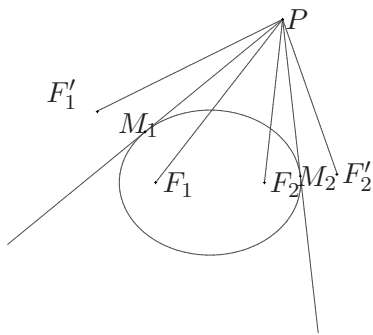,width=5cm}
\end{center}
\caption{The second little Poncelet theorem.}  \label{petitponc}
\end{figure}
\Pf
Suppose that $E$ is an ellipse (the case of a hyperbola is similar).

Reflect $F_1$ in $PM_1$ to $F'_1,$ and $F_2$ in $PM_2$ to $F'_2.$
As $ PM_1$ is a bisector of $ \widehat{F_1 M_1 F_2},$ we see that
$ F'_1, M_1, F_2$ are collinear and $F'_1 F_2= M_1F_1 +M_1 F_2 $ is the major axis of $E.$
We deduce that $F'_1F_2= F_1F'_2.$
Consequently, the triangles $F'_1PF_2$ and $F_1PF'_2$ are congruent, because
their sides are of equal length.
This implies $ \widehat{F'_1P F_2}= \widehat{F_1 P F'_2},$
and then
$$ \widehat{F'_1 P F_1} = \widehat{F'_1P F_2} - \widehat{F_1 P F_2}=
\widehat{F'_2 P F_1} - \widehat{F_1 P F_2 } = \widehat{F'_2 P F_2},
$$
which concludes the proof.
\EPf

\begin{theorem}
Suppose that  some segment of a billiard trajectory in an ellipse  does not intersect the
focal segment $ [ F_1 F_2].$
Then the billiard trajectory remains forever tangent to a fixed confocal ellipse  called
the caustic.
\end{theorem}
\Pf
\psfrag{A}{$P_1$}
\psfrag{B}{$M_1$}
\psfrag{C}{$P_0$}
\psfrag{D}{$M_2$}
\begin{figure}[th]
\begin{center}
\epsfig{file=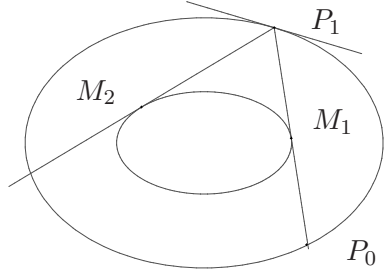,width=5cm}
\end{center}
\caption{Existence of an elliptic caustic.}  \label{caustic}
\end{figure}
Let $P_0 P_1$ be a segment of a billiard trajectory in an ellipse $E,$
and suppose that $ P_0 P_1$  does not intersect  $ [F_1 F_2].$

Reflect $ F_1$ in $P_0 P_1$  to $F'_1, $ and consider the ellipse
$ C= \{ MF_1 + MF_2= F_2 F'_1.$
We see that $M_1= F_2 F'_1 \bigcap P_0 P_1 $ belongs to $E,$
and since $ P_0 P_1$ is a bisector of $ \widehat{F_1 M_1 F_2},$  it is the tangent to $C$
at $M_1.$

Draw $P_1 M_2$ the second tangent to $ C.$
By the second little Poncelet theorem, the angles $\widehat{M_1 P_1 M_2}$
and $\widehat{ F_1 P_1 F_2}$ have the same bisectors.
Hence $P_1 M_2$ is the second segment of our billiard trajectory, and is tangent to $C.$

\EPf

\begin{remark}
When  some segment contains only one  focus, then every segment contains a focus, and there
is no caustic.
When some segment intersects the interior of the focal segment then there is a caustic,
which is a hyperbola with foci $ F_1$ and $ F_2.$
But the tangency points need not be at a finite distance. This fact will be illustrated
by the following example.
\end{remark}

\medskip\par\noindent
{\bf Example}\label{contrex}
\psfrag{P0}{$P_0$}
\psfrag{P1}{$P_1$}
\psfrag{P2}{$P_2$}
\psfrag{P3}{$P_3$}
\begin{figure}[th]
\begin{center}
\epsfig{file=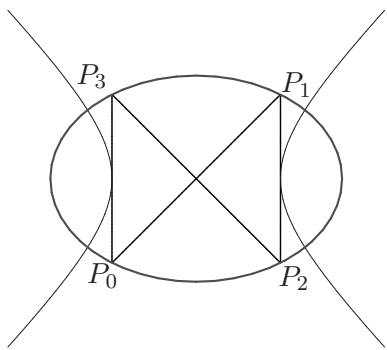,width=5cm}
\end{center}
\caption{A billiard quadrilateral in an ellipse. The caustic is a hyperbola.}  \label{contrex-dessin}
\end{figure}

Consider the points $ P_0= (-1,-1), \  P_1= (1,1), \ P_2 = (1,-1),$ and $P_3= (-1,1).$
The trajectory $ P_0P_1P_2P_3$ is a periodic billiard trajectory in the ellipse
$E= \{ x^2 + 2y^2 =3\}.$ The caustic is the hyperbola $ C=\{ x^2-y^2=1 \},$
and $ P_0 P_1, \  P_2 P_3$ are tangent to $C$ at infinity.

By Poncelet's theorem, if we apply the Poncelet construction   from another point $Q_0 \in E,$
then we obtain another
periodic billiard
trajectory $Q_0 Q_1 Q_2 Q_3.$  Moreover, the diagonals $Q_1Q_3$ and $Q_0 Q_2$ remain forever
parallel to the $x$-axis (this is a consequence of a theorem of Darboux that we will see later).
Consequently, all these Poncelet quadrilaterals are symmetric with respect to the
$y$-axis.

If we start from $Q_1=Q_3=A= (0 ,  \sqrt{3/2}), $ then the quadrilateral degenerates into
a trajectory of the form $ A\, B\, A\, C\, A,$
where $ B= ( \sqrt{ 5/3}, -\sqrt{2/3}), \   C= ( - \sqrt{ 5/3}, -\sqrt{2/3}).$

\medskip

The preceding example shows that a billiard trajectory $ A_0, \ldots , A_n$ such that $A_n=A_0$
needs not be $n$-periodic.
Moreover, a billiard trajectory in an ellipse is not necessarily a Poncelet polygonal
line. Consider for example an ellipse $E$ with foci $ F_1, F_2$ and a chord $AB$
such that $AB$ is the internal bisector of $\widehat{F_1 B F_2}.$
The polygonal line $ABA$ is a billiard trajectory in $E,$ but generally it is
neither a Poncelet polygonal line, nor a periodic billiard trajectory.

\smallskip

Happily, these defects cannot occur for billiard trajectories that do not intersect the
focal segment.

\begin{corollary}\label{periodic}
Let $P_0,P_1, \ldots , P_{n-1}, \, P_n= P_0$ be a billiard trajectory in an ellipse $E$ such that $P_0 P_1$
does not intersect the focal segment $ [F_1 F_2].$
Then it is a periodic billiard trajectory inscribed in $E$ and circumscribed about a confocal
ellipse $C.$
\end{corollary}
\Pf
Since $P_{n-1} P_0$ does not intersect the focal segment, its reflection at $P_0$ cannot be
$P_0 P_{n-1}.$ Since it is another tangent to $C$ through $P_0,$ it must be $P_0 P_1,$ and then
$P_{n+1}= P_1.$
\EPf

Now, we shall prove the existence of billiard polygons with $n$ sides and rotation number $p.$
 Following an idea due to Chasles and Birkhoff (\cite{Be,Bi}),
 we will obtain them as polygons of maximum perimeter among the $n$-polygons of
 rotation number $p$ inscribed in $E.$
We shall need the following classic lemma

\begin{lemma}\label{max}
Let $A,B$ be (not necessarily distinct) points of an ellipse $E.$ If the function
$ f(M)= MA + MB, \  M \in E$ has a local maximum at $C \in E,$ then
$CB$ is the reflection of $CA$ in the normal to $E$ at $C.$
\end{lemma}
\Pf
See \cite{LT}.
\EPf
\begin{proposition}\label{np}
Let $P_0$ be a point of an  ellipse $E.$ Let $n$ and $p$ be coprime integers such that
$n \ge 2p+1.$ There exists a billiard trajectory in $E,$ of period $n,$   winding number
$p$  and starting at $P_0.$
\end{proposition}
\Pf We shall consider the following maximum problem.
The domain of definition of the function to maximize is
$ {\cal A}= \{ (\alpha _1, \ldots , \alpha_n ) \in [0, \pi ]^n , \sum\alpha_i = 2 \pi p  \} ,$
it is a compact set.
For $ \alpha \in {\cal A}  $  and $P_0 \in E,$
let us define the inscribed polygon $ P_0, P_1, \ldots , P_n=P_0$ by the
angular condition
 $ \widehat{( \overrightarrow{FP_{i-1}}, \overrightarrow{FP_i})  }= \alpha_i,$
 where $F$ denotes a focus of $E.$
We want to maximize the perimeter  of this polygon, which is
$ f( \alpha)= P_0P_1 + P_1P_2+ \cdots + P_{n-1} P_0.$

Since ${\cal A}$ is compact and $f$ continuous, this maximum exists.
By lemma \ref{max}, it is a billiard trajectory at the points
$P_1, \ldots , P_{n-1}.$

Let us show that no segment $P_{i-1} P_i$ intersects the open focal segment
$(FF')$ of $E.$
If it was the case, then all segments of the trajectory would intersect $(F F'),$ hence
the winding number of this trajectory about $F$ would be zero, which is impossible by
the definition of our polygon.

\smallskip
If some segment was the major axis, then all segments would be this axis, and all
$ \alpha _i = \pi,$ which is impossible since $ \sum \alpha_i = 2p \pi < n \pi.$

\smallskip

If some segment contains only one  focus, then it is known that the billiard trajectory
converges to the major axis and is not periodic ( see \cite{St, Fr} ).

 We deduce
that $P_0 P_1$ does not intersect the focal segment $[F F']$ of $E.$
 By lemma \ref{max} and corollary \ref{periodic}, we conclude that $ P_0,P_1, \ldots P_{n-1} $
    is a periodic billiard trajectory
in $E.$

 The exact period $d$ of our trajectory
  is a divisor of $n,$ and we have $n=d \, u.$
By the angular  condition,  $u$ is a divisor of $p.$
Since $n$ and $p$ are coprime $u=1$ thus $n$ is the exact period of our trajectory.

\EPf
\begin{remark}
This does not  not prove that the caustics $C_n$ do not depend on the initial point $P_0.$
 This is true by Poncelet's theorem, which we shall prove later.
 Using a theorem of Graves {\rm \cite{Be} }, it can be shown
 that all the Poncelet polygons in two confocal ellipses have the same perimeter.
\end{remark}

\subsection{ Poncelet polygons and toric braids}

A toric braid is a braid  corresponding to the closed braid obtained by projecting
the standardly embedded torus knot into the $xy$-plane.
 A toric braid is a braid of the form
 $ \tau_{p , n} =  \Bigl(\sigma _1 \, \sigma _2 \, \cdots \sigma_{p-1} \Bigr)^n$,
where $ \sigma _1, \ldots, \sigma _{p-1}$ are the standard generators
of the full braid group $B_p.$

\begin{figure}[th]
\begin{center}
\begin{tabular}{ccc}
\epsfig{file=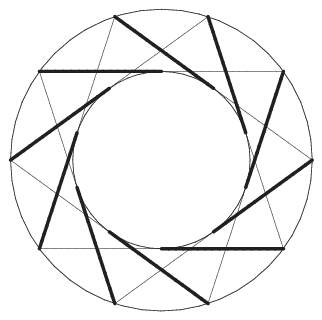,height=2.9cm,width=3.9cm}&
\epsfig{file=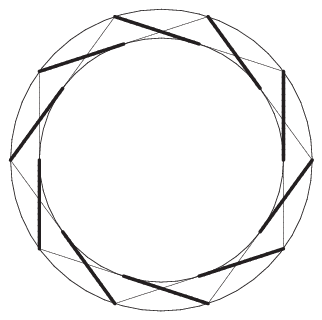,height=2.9cm,width=3.9cm}&
\epsfig{file=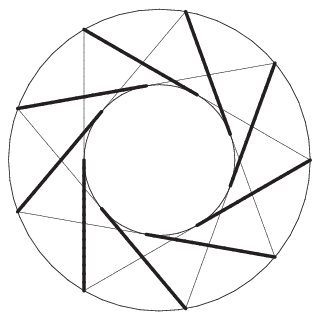,height=2.9cm,width=3.9cm}
\end{tabular}
\end{center}
\caption{Some Poncelet polygons (or unions of Poncelet polygons) in nested ellipses.
 They are projections of the toric braids
$\tau_{3,10}, \   \tau_{2, 10} $ and $\tau _{3,9} $,
and are denoted
 $ { 10 \brace 3}$, $ { 10 \brace 2}$  and $ { 9 \brace 3 }.$   }  \label{sp}
\end{figure}

\begin{remark}
Let $E$ and $C$ be nested ellipses such that there exists a Poncelet polygon
inscribed in $E$ and circumscribed about $C.$
 Every Poncelet polygon is the projection
 of a torus knot of type $ T( n,p), \  n\ge 2p+1.$ More precisely,
 if we cut  the elliptic annulus delimited by
 $E$ and $C$ along a half-tangent, then we see that such a polygon is ambient isotopic
  to the projection of the closure of the toric braid $ \tau _{p,n}.$
Consequently, it is also ambient isotopic  to the star
 polygon $ { n \brace p} ,$ see {\rm \cite{KP2}}.
\end{remark}

We shall need the following results on braids, due to Manturov \cite{M}.
A quasitoric braid of type $B (p,n)$
 is a braid obtained by changing some crossings in the toric braid
 $\tau_{p , n}.$

Manturov's theorem tells us that every knot (or link) is realized as the closure of
a quasitoric braid (\cite{M}).
 More precisely, he proved that any $\mu$-component link can be realized as the closure
 of a quasitoric braid of type $ B( p \mu, n \mu )$ where $ (p,n)=1,$  $p$ even and $n$ odd.

The quasitoric braids form a subgroup of the full braid group,  hence there exist trivial
quasitoric braids of arbitrarily great length. Hence
we can suppose $ n \ge 2p+1$ in Manturov's theorem.
Using this theorem, proposition \ref{np}, and Poncelet's theorem, we obtain
the main result of this section.

\begin{theorem}\label{projection}
Let $E$ be an ellipse. Every $\mu$-component link has a projection which is the union
of $\mu$   billiard trajectories in $E$ with the same odd period, and with the same caustic $C.$
\end{theorem}

\section{Jacobi's proof of Poncelet's theorem}\label{Jacobi}

We shall only need the following properties of elliptic functions, see \cite{WW} for
proofs.

\subsection{The Jacobian elliptic functions sn$\, {\bf z }$, cn$\, {\bf z}$ and dn$\, {\bf z}$.}

They depend on the choice of a parameter $k$, $0<k<1,$ called the elliptic modulus.

The Jacobi amplitude $\phi= {\rm am} (z) $ is defined by inverting the elliptic integral

$$ z= \int _0 ^{\phi} \frac{dt}{\sqrt{1-k^2 \sin ^2 t}}.$$
It verifies
$ {\rm am} (u+2nK)= {\rm am} (u) + n\pi,$ where $  {\rm am} (K)= \frac{\pi}{2},$ and
 $n \in {\bf Z}.$

The Jacobian elliptic functions are defined for $z$ real by
$$
{\rm sn} \, z=  \sin \bigl( {\rm am} (z) \bigr) , \quad
{\rm cn} \, z = \cos \bigl( {\rm am} (z) \bigr) , \quad
{\rm dn} \, z = \sqrt{1 - k^2 {\rm sn}^2 z},
$$
and can be extended to meromorphic functions on ${\bf C}.$
When $k=0,$ these functions degenerate into the
ordinary circular function $ \sin z$ and $\cos z.$
But, contrarily to the circular functions, they are doubly periodic functions with periods
$ 4K \in {\bf R}, $ and $4iK' \in i{\bf R},$
and they have poles.
For example, the poles of $ {\rm sn }\, z $ are congruent to
$iK' \   ({\rm mod.} {2K,2iK'}),$ its zeros are the points
 congruent to $0 \  ({\rm mod.}{2K,2iK'}),$
and its exact periods are $ 4 K, \ 2iK'.$
The zeros of $ {\rm cn}\, z$ are the points congruent to $ K  \  ({\rm mod.} {2K,2iK'}).$
We have
$ {\rm sn}(z + 2K) = - {\rm sn} \, z, \ $ and $ \ {\rm cn}(z + 2K) = - {\rm cn} \, z.$
 We also have ${\rm sn} (K+iK')= k^{-1},  $  which implies that
the zeros of ${\rm dn} \, z$ are the points congruent to $K+iK'  \ ({\rm mod.} {2K,2iK'}). $

We have the following addition formulas
$$
{\rm sn} (x+y) = \frac{{\rm sn} \, x \, {\rm cn} \,  y \, {\rm dn} \, y +
{\rm sn} \, y \, {\rm cn}\,  x\,  {\rm dn} \, x   }
{1- k^2 \, {\rm sn}^2 \, x \, {\rm sn} ^2 \,y}, \  {}  \
{\rm cn} (x+y)= \frac{ {\rm cn} \, x\, {\rm cn} \, y \, -
{\rm sn} \, x \, {\rm sn} \, y \, {\rm dn} \, x \, {\rm dn} \, y }
{1- k^2 \, {\rm sn}^2  x \, {\rm sn} ^2 y}
$$
When $k=0,$ these formulas degenerate into the usual addition
 formulas for the circular functions.

In the next section we will use   the following formula
due to Jacobi (\cite{WW} p.529).
$$
\sin \bigl( {\rm am} (u+v) + {\rm am} (u-v) \bigr) =
\frac{2 \, {\rm sn } \, u \, {\rm cn} \, u \, {\rm dn } \, v }
{1- k^2 \, {\rm sn}^2  u  \, \, {\rm sn} ^2  v}
$$

\subsection{ Jacobi's uniformisation }

The next result is a variant  of Jacobi's
 uniformization of the Poncelet problem. It is due to Hermite and Laurent (\cite{Lau}).

\psfrag{A}{\small $P(\varphi+\beta)$}
\psfrag{B}{\small $M(\varphi)$}
\psfrag{C}{\small $P(\varphi-\beta)$}
\psfrag{D}{\small $P(\varphi+2\beta)$}
\begin{figure}[th]
\begin{center}
\epsfig{file=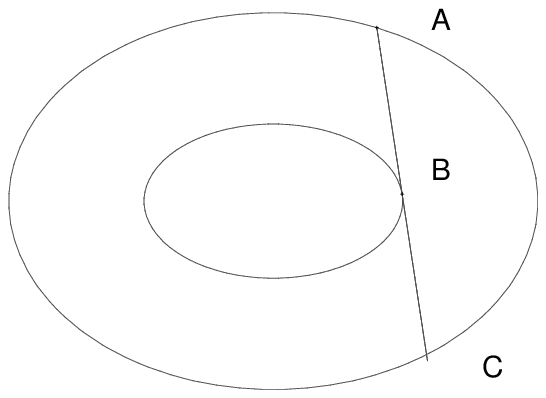}
\end{center}
\caption{Jacobi's lemma.}  \label{jacobi}
\end{figure}

\begin{lemma}
Let $E$ and $C$ be the ellipses defined by
$$ E= \Bigl\{ \frac{x^2}{a^2} + \frac{y^2}{b^2}=1 \Bigr\} , \ a>b>1, \ {} \  {}  \
 C= \Bigl\{ x^2+y^2=1 \Bigr\} .$$
Let us parameterize $E$ by $P( \psi ) = ( a \, {\rm cn} \, \psi , \, b \, {\rm sn} \, \psi ),$
and $C$ by $M( \phi) = ( {\rm cn} \, \phi , \, {\rm sn} \, \phi ),$
where the elliptic modulus $k$ is defined by $ k^2 ( a^2-1) = (a^2-b^2).$
Let $\beta$ be a real number such that $ {\rm cn} \, \beta = 1/a.$

Then the tangent to $C$ at $ M(\phi)$ intersects $E$ at
 $ P( \phi - \beta)$ and $ P( \phi + \beta).$
\end{lemma}
\Pf
We have $ {\rm dn}^2  \beta= 1 - k^2 \, {\rm sn}^2   \beta = b^2/a^2,$
hence $ {\rm dn} \, \beta =b/a.$

Let us show that $ P( \phi + \beta)$ belongs to the tangent to $C$ at $M( \phi).$
The equation of this tangent is
$ x \, {\rm cn } \, \phi + y \, {\rm sn} \, \phi =1. \ $
Let us compute
$S = a \, {\rm cn}(\phi + \beta) \, {\rm cn} \, \phi +
 b \, {\rm sn}( \phi +\beta)\, {\rm sn}\, \phi.$

Using the addition formulas we obtain
$$
S \, ( 1- k^2 \, {\rm sn}^2  \phi \, {\rm sn}^2  \beta ) =
{\rm cn}^2 \phi + {\rm sn}^2  \phi \, {\rm dn}^2  \beta =
{\rm cn}^2  \phi + {\rm sn}^2  \phi \, ( 1 - k^2 \, {\rm sn}^2 \beta )=
1-k^2 \, {\rm sn}^2  \phi  \, {\rm sn}^2  \beta.
$$
Consequently, $S=1,$ and $P( \phi + \beta )$ belongs to the tangent to $C$ at $M( \phi).$
Changing $\beta$ to $- \beta,$ we see that $P (\phi - \beta)$ also belongs to this tangent.
\EPf

\begin{remark}
By affinity, Jacobi's lemma extends easily in the
case of two nested ellipses with the same two axis, meeting transversally in $P_2 ( {\bf C}) .$
 When this pair of ellipses
is affinely equivalent to a pair of concentric circles, the elliptic parametrizations degenerate
into the usual circular ones.
\end{remark}

\subsection{Proof of Poncelet's closure theorem}

\begin{figure}[th]
\begin{center}
\epsfig{file=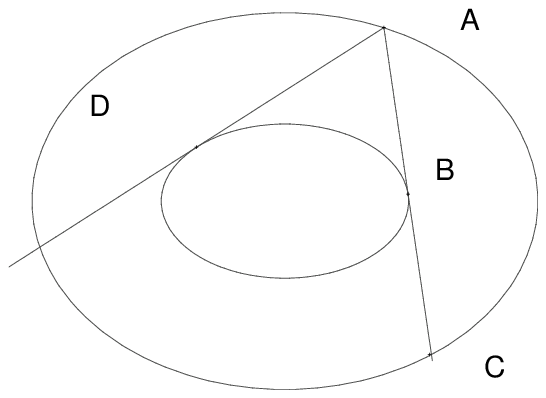}
\end{center}
\caption{Proof of Poncelet's closure theorem.}  \label{demponc}
\end{figure}

We shall now present  the  Hermite--Laurent proof of Poncelet's theorem for
a pair of confocal ellipses. Since any pair of conics meeting  transversally
 in $ P_2( {\bf C})$
is projectively equivalent to a pair of confocal ellipses ( \cite{LT}),
  we obtain a proof of the generic case of Poncelet's theorem.
For the nongeneric cases see \cite{Be,Sa}, and for
the original proof of Jacobi see \cite{BKOR}.
\Pf
Let $P_0, P_1, \ldots P_{n-1} P_0$ be a Poncelet polygon inscribed in $E$ and circumscribed
about $C.$
Let $P_j P_{j+1} $ be tangent to $C$ at $M_j.$
We will use the Jacobi parametrizations of $E$ and $C.$
If $M_0= M( \phi), $ then by Jacobi's lemma we can suppose $P_1= P( \phi  + \beta).$
Using Jacobi's lemma again, we have  $M_1= M( \phi + 2 \beta),$
and by induction  $M_j= M ( \phi + 2j \beta ).$

Since the polygon closes after $n$ steps, we have $M_n= M_0,$ or
$  M( \phi + 2n \beta )=M( \phi ) .$

That means
$ {\rm am} ( \phi + 2n \beta) = {\rm am} \, \phi + 2q \pi = {\rm am} ( \phi + 4 q K)$
 by the properties of the Jacobi amplitude.
Consequently we obtain $ 2n \beta = 4 q K,$ or  $\beta = 2q K /n.$

Now, let us consider a Poncelet polygon starting from an arbitrary point
$ M'_0 = M( \phi ' ) $ of $C.$
By Jacobi's lemma we have $ M'_n= M( \phi ' +2n \beta) = M(\phi ' +4q K ) = M( \phi ') = M'_0.$

Consequently, we see that every Poncelet polygonal line closes after $n$ steps.
\EPf

\begin{corollary} Let $ E$ and $C$ be confocal ellipses, and let ${\cal P}$
be a Poncelet polygon with an even number of sides. Then $ {\cal P}$ possesses
a central symmetry.
\end{corollary}
\Pf
We use the same notation as before. We have $ \beta= 2q K/n$ ,
where  $n=2h$ is the number of sides of $ {\cal P}. $
The numbers $n$ and $q$ are coprime, hence $q$ is odd.
For every $\phi,$ we have
$$ M_h (\phi) = M ( \phi + 2h \beta) = M(\phi + 2q K)= M( \phi + 2K)=
  \bigl(  {\rm cn} (\phi + 2K),    \  {\rm sn} ( \phi + 2K ) \bigr)= - M_0 (\phi).$$
\EPf

By projective equivalence,
this implies a remarkable theorem of Darboux    (\cite{Da}).

\begin{theorem}
Let $ {\cal P}$ be a Poncelet polygon with an even number of sides inscribed in a conic $E$
and circumscribed about another conic $ C$.
Then the diagonals of ${\cal P}$ pass through a point,
which is the same for every  Poncelet polygon.
\end{theorem}


\section{Irregularity of Poncelet odd polygons}

Most regularity properties of a polygon can be expressed by rational linear relations
 between some of its segments.
Let us parameterize a (crossed) polygon by arc length, starting at a vertex $P_0.$
 We shall say that this polygon is completely irregular if $1$ and the arc lengths of
its crossings and  vertices (except $P_0$) are linearly independent over $ {\bf Q}.$

The purpose of this section is to prove that if $E$ and $C$ is a pair of confocal ellipses
possessing a Poncelet polygon with an odd number of sides, then there exists a completely
irregular Poncelet  polygon.  We will give an analogous result for unions of finitely
many Poncelet polygons.

\subsection{ Two lemmas on elliptic functions}
We shall use elliptic functions to compute the arc lengths of the crossings and vertices
of Poncelet polygons. We shall need the following two technical lemmas.

\begin{lemma} \label{lem1}
Let $n$ and $p$ be coprime integers, with $n$ odd.
For every integer $j,$ let us define the function
$ f_j (z) = {\rm sn}^2 (z + j \theta) + r^2,$ where $r^2>0$  and $ \theta= 4 p K / n.$
Then, if $ h \not\equiv j    \  ( \mod n ),$  the functions
$f_j (z)$ and $f_h (z)$ do not possess any common zero.
\end{lemma}
\Pf
First, let us study the zeros of the elliptic function
$ g(z) = {\rm sn} \, z + i r, \  r>0.$ By considering its restriction to the $y$-axis,
we see that there exists a pure imaginary $ \alpha, $ such that $ g(\alpha)=0.$
Since we have $ {\rm sn}  (2K - \alpha) = {\rm sn} \,  \alpha ,$ we see that
$2K- \alpha$ is another zero of $g(z).$  As $g(z)$  is an elliptic function of order two,
its  zeros are the points congruent to
$ \alpha$  or $2K - \alpha   \ ( {\rm mod .} \ 4K, 2iK').$
By parity, we deduce that the zeros of $ f_j(z)$ are the numbers which are congruent
to $ \pm \alpha - j \theta, $ or $ 2K \pm \alpha - j \theta, \  ( {\rm mod .} \ 4K, 2i K' ).$

If we had $ \alpha - j \theta  \equiv \alpha - h \theta,  $  or
$ \alpha - j \theta \equiv  2K + \alpha - h \theta \  ( {\rm mod .}  4K, \, 2iK'),$
then we would deduce
$ (h-j) \theta \equiv 0  \ ( \mod 2 K) .$
This implies that $ 2(h-j)p/n$ is an integer, which is impossible since $n$ is odd,
$ (n,p)=1,$ and $h  \not\equiv j \ ( \mod n) .$

If we had $ \alpha - j \theta \equiv - \alpha - h \theta$ or
$ \alpha - j \theta \equiv 2 K -  \alpha - h \theta, \  ( {\rm mod .}  4K, \, 2iK'),$
then we would have
$ 2 \alpha \equiv  ((j-h) \theta \  ( {\rm mod .} \ 2K, 2iK' ).$
Taking the real parts, we would obtain
$ (j-h) \theta \equiv 0 \ ( \mod 2K )$ which is impossible.

Consequently, $\alpha - j \theta $ cannot be a zero of $f_h (z).$
The proof that the other zeros of $f_j (z) $ cannot be zeros of $f_h(z)$ is entirely similar.
\EPf

\begin{remark}
As the proof shows it, the condition $n$ odd is  necessary in  lemma \ref{lem1}.
\end{remark}

\begin{lemma} \label{lem2}

Let $n$ and $p$ be coprime integers.
For $ j \not \equiv 0 \ (  \mod n  ) ,$ let us define the functions $D_j (z)$ and $F_j (z)$ by
$$
D_j (z)= {\rm sn}(z+j \theta)\, {\rm cn} \, z  -{\rm cn}(z+j \theta)\, {\rm sn} \, z , \  {} \
F_j(z)= \frac{ {\rm sn}( z + j \theta ) - {\rm sn} \, z }{D_j (z)}, \  {\it where} \
\theta = \frac{4 p K}{n}.
$$
 Then, for every integer $j$ there exists a complex number $ \alpha_j$ such that
 $ F_j ( \alpha_j)=\infty, $ and if $h \not\equiv j \ ( \mod n ) ,$
 then $  \   F_h (\alpha _j) \ne \infty.$
\end{lemma}
\Pf
We have
$$ D_j (z) = \sin \bigl( {\rm am} ( z + j \theta) - {\rm am } \, z \bigr)=
\sin \bigl({\rm am} ( z + j \theta) + {\rm am }(- z) \bigr)
$$
Now, using the Jacobi formula for $ \sin \bigl( {\rm am }(u+v) + {\rm am} (u-v) \bigr),$
we obtain
$$ D_j (z)=\frac{ 2 \, {\rm sn} (j \beta) \, {\rm cn}(j \beta)  \,
{\rm dn} (z + j \beta)}
{1 - k^2 \, {\rm sn}^2 ( j \beta)  \, {\rm sn}^2 (z+ j \beta ) },  \  {} \ {}
\  {\rm where} \ \beta= \frac{\theta}{2}
$$
Let $ \alpha_j = -j \beta + K + iK' .$
We have $ {\rm dn} ( \alpha _j + j \beta)= {\rm dn} ( K+iK')=0.$
Since $ {\rm dn}^2 z + k^2 \, {\rm sn}^2 z =1, $ we obtain
$ {\rm sn}^2( \alpha_j + j \beta) = 1/k^2, $ and then $D_j (\alpha_j )=0.$

The numerator of $F_j ( \alpha_j) $ is
$$N( \alpha _j) = {\rm sn} ( \alpha _j + j \theta ) -{\rm sn} ( \alpha _j )=
{\rm sn} ( K+iK' + j \beta ) - {\rm sn} ( K+iK' - j \beta ).$$

Using the addition formula for the function $ {\rm sn} \, z,$ we obtain
$$ N( \alpha _j) = 2 \frac{ {\rm sn} (K+iK') \, {\rm cn}( j \beta) \, {\rm dn} (j \beta)}
{1 - k^2 \, {\rm sn}^2(K+iK') \, {\rm sn}^2( j \beta)} .
$$
Since $ {\rm sn} (K+iK')= k^{-1},$ we obtain
$ N (\alpha _j) = 2 k^{-1} \Frac{  {\rm dn}( j \beta)}  { {\rm cn} ( j \beta) } \ne 0,$
and then $ F_j ( \alpha _j) = \infty.$

\medskip

On the other hand, if $ h \not\equiv j \  (\mod n) ,$ we have
$ \alpha _j + h \beta= K+ i K' + 2(h-j)p K /n.$

First, we see that $ \alpha _j + h \beta \not\equiv  K+iK' \ ({\rm mod .} 2K, 2iK' ) ,$
which implies that
$ {\rm dn } ( \alpha _j + h \beta ) \ne 0.$

We also see that $ \alpha _j + h \beta \not\equiv  iK' \ ({\rm mod .} 2K, 2iK' ),$
which implies that $ {\rm sn } ( \alpha _j + h \beta ) \ne \infty .$
We conclude that $ D_h ( \alpha _j) \ne 0.$

\smallskip

Let us show that if $ {\rm sn } \, z =\infty,$ then $ F_h (z) \ne \infty .$

Since the functions $ {\rm sn}\, z$ and $ {\rm sn} (z + h \theta)$ do not have common poles,
 $ {\rm sn} (z + h \theta) \ne \infty .$

On the other hand, as $ {\rm sn}^2 z + {\rm cn}^2 z =1,$
we obtain
$$
\frac{ {\rm cn}^2 \,  z  }{  {\rm sn} ^2 \,z}=-1,  \  {\rm and \  then} \
F_h(z)= \frac{-1}
{  {\rm sn} (z+ h \theta) \Frac{ {\rm cn } \, z }{ {\rm sn } \, z} - {\rm cn} ( z+ h \theta)}
$$
If we had $ F_h(z)= \infty,$ then
$$
{\rm sn ( z + h \theta) \frac{ {\rm cn } \, z }{ {\rm sn } \, z}
 ={\rm cn} ( z+ h \theta)},
$$
whence
$ {\rm sn}^2 (z +h \theta)= - {\rm cn}^2 (z +h \theta) \ne \infty,$
and $ {\rm sn}^2 (z +h \theta)+ {\rm cn}^2 (z +h \theta)= 0,$ which is impossible.

Similarly, we see that if
$ {\rm sn} (z +h \theta )=\infty, $ then $F_h(z)  \ne \infty.$

Now, let us prove that $ F_h (\alpha _j ) \ne \infty.$
We have $ D_h ( \alpha _j) \ne 0,$ and we have  proved that we can suppose
$ {\rm sn} (\alpha_j ) \ne \infty $ and $ {\rm sn} (\alpha _j  +h \theta )   \ne \infty, $
then
$$
F_h ( \alpha _j) = \frac{ {\rm sn} (\alpha _j + h \theta) - {\rm sn} \, \alpha _j }
{ D_h ( \alpha _j )} \ne \infty.
$$
\EPf
\subsection{  Irregular Poncelet polygons with an odd number of sides}
\begin{proposition}\label{arclengths}
Let $E$ and $C$ be confocal ellipses such that there exists a Poncelet polygon  $ {\cal P}$
 inscribed in $E$ and circumscribed about $C.$
We suppose that the number of sides of $ {\cal P}$  is odd.
Then there exists a Poncelet polygon  satisfying the following condition.

If the arc lengths $t_i$ of the vertices and crossings are measured from a vertex $P_0,$
then the numbers $1$ and $t_i, \ t_i \ne 0$ are linearly independent over ${\bf Q}$.
\end{proposition}
\Pf
$$
{\rm Let} \   E = \{ \frac{x^2}{a^2} + \frac{y^2}{b^2}=1 \} , \  a > b >1 , \ {\rm and} \
\ C = \{ x^2 + \frac{y^2}{c^2} = 1 \} , \  c<1
$$
be our ellipses. The condition on the eccentricity of $ C$ means that $2 c^2 >1.$

Let us consider the Jacobi parametrizations of $E$ and $C$ by means of elliptic functions,
and let $ \theta = 4 p K / n .$
To each real number $\phi$ corresponds a Poncelet polygon $ {\cal P} _ {\phi} $
through $ M( \phi) = \bigl( {\rm cn } \, \phi , \,  c \, {\rm sn} \, \phi \bigr) .$
Let us denote $ \phi _j = \phi + j \theta, \  M_j= M(  \phi _j ), $ and let
$ \ell_j$ be the tangent to  $ C$  at $M_j.$
The equation of $ \ell_j$ is
 $$ x \, {\rm cn} \, \phi _j + \frac{y }{c } \, {\rm sn} \, \phi _j = 1.$$

Let $ Q _{h,j}= \ell _h \bigcap \ell _{h+j} , \  j \not\equiv 0 \  ( \mod n ).$
The abcissa $x_{h,j}$ of $Q _{h,j}$ is
$$
x_{h,j} = \frac
{- \, {\rm sn} \, \phi _h +  {\rm sn} ( \phi_h + j \theta) } { {\rm sn} (\phi _h + j \theta)
\, {\rm cn} \, \phi _h - {\rm cn} (\phi _h + j \theta) {\rm sn} \, \phi_h  }
= F_j (\phi _h)
$$
where $F_j$ is the function defined in lemma \ref{lem2}.
The abcissa of  $ P_h = Q_{h, -1} = Q _{ h-1, 1}$ will also be denoted by $ x_h= x_{h,-1}.$
The distance $ P_h Q_{h,j}$ is $ | d_{h,j}|$ where
$$
d_{h,j}= d_{h,j} ( \phi) = \frac{\sqrt{1-c^2}}{ {\rm sn}\, \phi_h} \
\sqrt{ {\rm sn}^2 \phi _h + \frac{c^2}{1-c^2}} \
\Bigl( x_h - x_{h,j} \Bigr)
$$
Since $ c^2 / (1-c^2) >1 ,$ the function $d_{h,j} ( \phi) $ is meromorphic in a neighborhood
of the real axis.

\medskip

Our first step is to prove that the functions $1$ and
$ d_{h,j} ( \phi), \ j \not\equiv -1 \ ( \mod n ) $ are linearly independent over $ {\bf C}.$

Let $\lambda _{h,j} $ and $ \lambda $ be complex numbers such that
$ \sum _{h=1}^{n } \, \sum _{j=1} ^{n-2} \, \lambda _{h,j} d_{h,j} = \lambda,\ $ or
$$
\sum _{h=1} ^n \frac{\sqrt{1-c^2}}{ {\rm sn} \, \phi _ h } \
\sqrt{ {\rm sn}^2 \phi _h + \frac{c^2}{1-c^2}} \
\Biggl(  \, \sum_{j=1} ^{n-2} \lambda _{h,j} (x_h- x_{h,j} ) \Biggr) = \lambda.
$$

Since  $ c^2 /(1-c^2) >0, $  we see by lemma \ref{lem1} that the functions
$ f_h ( \phi)  = \sqrt{  {\rm sn }^2  \phi_h + c^2/(1-c^2)  }$ do not possess any common zero.
Hence, in the neighborhood of a zero of $ f_h ( \phi)$
this function is not meromorphic, while the other functions  are.

This implies that for every $h = 1 \ldots n $ we have
$$
\sum _{j=1} ^{n-2} \lambda _{h,j} \bigl( x_h- x_{h,j} \bigr) =0, \  {\rm and \  then} \
\lambda=0.
$$
Using our expressions of the abcissas $x_{h,j},$ we obtain the following
 relation between meromorphic functions
$$
\sum _{j=1} ^{n-2} \lambda _{h,j} \bigl( F_{-1} (z) - F _j (z) \bigr) =0.
$$
By lemma \ref{lem2},  for every integer $j \ne 0$ there exists a number $ \alpha _j$
such that $ F_j ( \alpha _j) = \infty, $ and $ F_h ( \alpha _j ) \ne \infty  \  $
if $ h \not\equiv j \ ( \mod n ) .$
Letting $ z= \alpha_j,$ we obtain $ \lambda _{h,j}=0,$ which concludes the proof of the
linear independence  of our functions.

\medskip

Now, we shall prove that for most $ \phi \in {\bf R},$ the numbers
$ d_{h,j} ( \phi)$ and $1$ are linearly independent over $ \QQ.$

For every nonzero collection of rational numbers $  \Lambda = ( \lambda, \lambda _{h,j}),$
let us define the function $ F_{ \Lambda} $ by
$ F_{ \Lambda} ( \phi) = \lambda - \sum _{h,j} \lambda _{h,j} \, d _{h,j} ( \phi).$
By our first step, this function is not identically zero, and it is meromorphic in a
neighborhood of ${\bf  R}.$  Therefore, the set of its real zeros is countable.
Consequently, the set of all real numbers $ \phi$ such that
$1$ and the numbers $d_{h,j} ( \phi) $ are linearly dependent over $ {\bf Q} $
is countable.
By cardinality, we deduce that the complementary set is not  countable, hence nonempty.
Consequently, there exists a real $ \phi$ such that $1$ and the numbers
$ | d_{h,j} ( \phi) |$ are linearly independent over $ {\bf Q}.$

\smallskip

Now, let us parameterize our Poncelet polygon by arc length, starting from $P_0$
for $t_0 \in {\bf Q} .$
The arc length $t_{h,j}$ of $ Q_{h,j} $ is
$$\displaylines {
t_{h,j} = t_0 +  d(P_0, P_1) +  d (P_1, P_2) + \ldots +  d( P_{h-1}, P_h)
+ d ( P_h, Q_{h,j} )\cr
= t_0 + | d_{0,1}| + |d_{1,1}| + |d_{2,1}| + \ldots + |d_{h-1,1}| +
|d_{h,j}|.\cr}
$$

The result follows from the independence of the numbers $1$ and $ |d_{h,j}|.$
\EPf

We shall also need an analogous result for links.

\begin{proposition} \label{links}
Let $E$ and $C$ be confocal ellipses such that there exists
a  polygon of an odd number of sides inscribed in $E$ and circumscribed about
$C.$

For any integer $\mu$,  there exist $\mu$ Poncelet polygons
$ {\cal P}^{(0)}, {\cal P}^{(1)}, \ldots , {\cal P}^{( \mu-1)} $
 satisfying the following condition:

for each such polygon, if $t_i$ are the arc lengths corresponding to its vertices , its
crossings, an its intersections with the other polygons, then the numbers
$ 1$ and $t_i, \ i\ne 0$ are linearly independent
over $\QQ.$
\end{proposition}
\Pf
Let $ \tau = \Frac {\theta}{\mu}, $ and let us denote
 $M_h= M ( \phi + h \tau) \in C,$ and $\ell_h $ the tangent to $C$ at $M_h.$
Let us consider the Poncelet polygons
$ {\cal P}= {\cal P}^{(0)}, \  {\cal P} ^{(1)},  \ldots , {\cal P}^{( \mu-1)}$ through
the points $ M_0, M_1, \ldots , M_{\mu-1}.$
The polygon $ {\cal P}$ is tangent to $C$ at the points
$ M_0, M_{\mu}, M_{2 \mu} , \ldots , M _{(n-1) \mu}.$
The vertices and crossings of  ${\cal P}$ are the points
$ Q_{h,j} = \ell _h \bigcap \ell_{h+j} , $ where $ h \equiv 0 \ ( \mod  \mu ) .$

Just as before, it can be proved that the distances $1$ and
$ | d_{h,j} ( \phi ) |, \   h \equiv 0,  \  j \not\equiv 0 , \  j\not\equiv -1  \  ( \mod \mu )$
are linearly independent over $\QQ,$ except for a countable set of numbers $ \phi.$

Consequently, the number $1$ and the arc lengths $t_i, \ i \ne 0 $ of
 the crossings and vertices of $ {\cal P}$ are linearly independent over
 $ \QQ$ except on a countable set of values of $\phi.$

By cardinality, we can suppose that the same property is true for each polygon
$ {\cal P}^{(j)}, \ j=0, \ldots ,  \mu -1,$ which proves our result.
\EPf
\section{Proof of the theorem }
We will  use Kronecker's theorem ( see \cite[Theorem~443]{HW}): 
\begin{theorem}
If $ \theta _1, \theta_2, \ldots, \theta_k, 1$ are linearly independent over $\QQ,$
then the set of points
$ \Bigl((n \theta _1), \ldots , (n \theta_k) \Bigr) $ is dense in the unit cube.
Here $ (x)$ denotes the fractional part of $x.$
\end{theorem}
Now, we can prove our main theorem.
\begin{theorem}\label{th:billiardkl}
 Let $E$ be an ellipse which is not a circle,
 and let $D$ be the elliptic cylinder $ D= E \times [0,1].$
Every knot (or link) is a billiard knot (or link)  in $D$.
\end{theorem}
\Pf
First, we consider knots.
By theorem \ref{projection} there exists a knot  isotopic to $K,$
whose projection on the $xy$-plane is a billiard trajectory of odd period
in the ellipse $E.$
If $t_0, t_1, \ldots , t_k $ are the arc lengths corresponding
to the vertices and crossings, we can
suppose by proposition \ref{arclengths} that the numbers $ t_1, \ldots , t_k,$
 and $ 1$ are linearly independent over $ \QQ$.
Using a dilatation, we can suppose that the total length of the trajectory is $1.$

Let us consider the  polygonal curve  defined by $ (x(t), y(t),z(t)),$
where $z(t)$ is the sawtooth function $z(t)= 2 | (nt+ \phi)-1/2 |$
 depending on the integer $n$ and on the real number $\phi.$
If the heights $z(P_j)$ of the vertices are such that $z(P_j) \ne 0,  \ z( P_j) \ne 1$, then
it is a periodic billiard trajectory in the elliptic cylinder $  {\bf D} =  E \times [0,1] $
(see \cite{JP,La2,LO,P,KP1}).
If we set $ \phi = 1/2 + z_0/2,$ $ z_0 \in (0,1),$ we have $ z(0)=z_0.$
Now, using Kronecker's theorem, there exists an integer $n$ such that
the numbers $z(t_i)$ are arbitrarily close to any chosen collection of heights,
which completes our proof.
\pn
The case of $ \mu$-component links is similar. First, by theorem \ref{projection},
we  find a  diagram that   is the union of $\mu$ Poncelet polygons
 with the same odd number of sides.
  Then, by proposition \ref{links} and Kronecker's theorem,
 we  parameterize each component so that the heights of the vertices and crossings are close
to any  chosen numbers.
\EPf

\vfill
\hrule
\pn
D. Pecker, U. Pierre et Marie Curie (Paris 6), Math\'ematiques, {\tt pecker@math.jussieu.fr}
\end{document}